\documentclass[preprint,12pt,authoryear]{elsarticle}

\usepackage{amssymb}
\usepackage[utf8]{inputenc}
\usepackage[english]{babel}
\usepackage{amsmath}
 \usepackage{graphicx}

\newtheorem{property}{Property}

\journal{European Journal of Operations Research}

\begin{document}

\begin{frontmatter}



\title{A graph-based algorithm for the multi-objective optimization of gene regulatory networks}


\author{Philippe Nghe*}
\address{FOM Institute AMOLF, Science Park 104, 1098XG Amsterdam, the Netherlands}
\address{Laboratoire de Biochimie, UMR CNRS-ESPCI Chimie Biologie Innovation 8231, PSL Research University,  ESPCI, 10 rue Vauquelin, 75005 Paris, France}
\address{*corresponding author}

\author{Bela M. Mulder, Sander J. Tans}
\address{FOM Institute AMOLF, Science Park 104, 1098XG Amsterdam, the Netherlands}

\begin{abstract}
The evolution of gene regulatory networks in variable environments poses Multi-objective Optimization Problem (MOP), where the expression levels of genes must be tuned to meet the demands of each environment. When formalized in the context of monotone systems, this problem falls into a sub-class of linear MOPs. Here, the constraints are partial orders and the objectives consist of either the minimization or maximization of single variables, but their number can be very large. To efficiently and exhaustively find Pareto optimal solutions, we introduce a mapping between coloured Hasse diagrams and polytopes associated with an ideal point. A dynamic program based on edge contractions yields an exact closed-form description of the Pareto optimal set, in polynomial time of the number of objectives relative to the number of faces of the Pareto front. We additionnally discuss the special case of series-parallel graphs with monochromatic connected components of bounded size, for which the running time and the representation of solutions can in principle be linear in the number of objectives.

\end{abstract}


\end{frontmatter}


\section{Introduction}

Adaptation to changing environments \textit{via} gene regulation is by nature a multi-objective problem, where the expression level of genes must be optimally set given multiple possible combinations of environmental signals \citep{Poelwijk2011}. The multi-valued output response to a collection of inputs is determined by the connectivity of the gene network  \citep{Alon2006}. Biological evolution adds a second layer of complexity, whereby the strength of the connections can be altered by mutations, modulating the multi-valued response itself. There does not currently exist any generic approach as to predict bounds to how well this response can be optimized, given the connectivity of the network.

Here, we apply the notion of Pareto optimality to gene expression in several environments. Pareto optimality is a natural extension of the concept of maximum to
multi-objective optimization problems. A solution is part of the Pareto
optimal set, or Pareto front, if it is impossible to improve one objective
without worsening another. Instead of imposing an aggregation of the
different objectives into a scalar function, Pareto optimality keeps
track of all potentially interesting solutions in the presence of trade-offs.
The Pareto approach, originally introduced in economics
\citep{Pareto1906,Voorneveld2003}, has proved useful in many engineering
applications \citep{Ehrgott2000, Zitzler2003, Geilen2005}, decision-making
analysis \citep{Yang2003}, and recently, medecine \citep{Cruz2013} and biology \citep{Alon2012, Schuetz2012}.

We approach the problem in the framework of monotone systems, which is widely used in control theory \citep{Angeli2003}, and more specifically for modelling gene regulatory networks \citep{Sontag2005}. This formalism, as detailed in section 2, leads us to define partial order constraints between the expression levels of a gene given the environmental inputs. Expression of a gene is considered to be either beneficial or detrimental in each given environment, meaning that expression of a gene must be either minimized or maximized given an input. The expression level of the gene in each environment corresponds to a dimension of the decision space. Given that one objective is associated with each environment, the objective space has the same dimensionality.

The partial order constraints that define the feasible set imply that our problem falls into the category of linear Multi-objective Optimization Problems (linear MOPs) \citep{Greco2005}. This problem could in principle be tackled with existing strategies, such as the multi-objective simplex algorithm or Benson's algorithm \citep{Ehrgott2012, Lohne2017}. However, although such algorithms can cope with large sized decision spaces, in our case the number of optimization objectives equals the number of decision variables and can reach several dozen. The unusually high number of objectives imposed by our application is considered a particularly difficult problem in the general case, as the category of \textit{many-objectives optimization problems} starts as soon as there are more than three or four objectives \citep{Fleming2002,Jaimes2015}. 

The algorithm presented here exploits the specificities of our linear MOP, in which constraints are exclusively partial orders and objectives are either the maximization or minimization of coordinates of the decision variables. This algorithm provides an efficient, exact and exhaustive description of the Pareto front, even with such a large number of objectives.

Note that the resolution of the problem provided in the present work could apply beyond the framework of monotone systems, in cases where constraints are expressed in the form of partial orders from the start. This could be the case for task scheduling problems: some tasks must be realized in a certain temporal order relative to each other due to design constraints \citep{Policella2007} (e.g. the assembly of the different parts of a car), with some having to be carried out as soon as possible and others as late as possible due to externalities (e.g.\ supply constraints, processing unit occupancy). 

This paper is structured as follows: in section 2, we detail the problem of determining bounds to the evolutionary potential of gene regulatory networks and its formalization in terms of multi-objective optimization under partial order constraints. In section 3, we provide a graph formulation of this problem, using Hasse diagrams \citep{Skiena1990}, which we colour according to optimization objectives. In section 4, we describe and prove a graph algorithm based on successive edge contractions with appropriate vertex colouring rules. In section 5, we discuss the complexity of our algorithm as a function of the number of objectives $N$ and of the number $N_P$ of maximal convex subsets of the Pareto front, or maximal efficient faces as defined in the papers \citep{Yu1975} or \citep{Ecker1980}. This leads us to propose an improved version of the algorithm running in O($P(N) \cdot N_P)$ time, where $P$ is a polynomial. We discuss how a parameterized complexity approach \citep{Alber2004} can provide a combinatorial description of the Pareto front, with a complexity of $O(N)$ in the case of series-parallel partial orders with monochromatic connected components of bounded size. In section 6, we provide an example of an explicit resolution and an exhaustive characterization of Pareto front sizes in the problem of 2D dimensional gradients of signals with stress patches.

\section{Multi-objective optimization in regulatory networks}

Organisms are typically confronted with a large variety of environmental signals that can themselves combine into an even larger diversity of spatio-temporal niches \citep{Chait2016}. Whether organisms are able to evolve an appropriate response to this diversity of environments depends fundamentally on the constraints imposed by their current regulatory responses. Despite the multi-objective nature of this problem, research on the evolution of gene regulatory networks has so far relied almost exclusively on explicit simulations of complex responses given a single objective, such as the number of maxima of a biological trait in a gradient \citep{Francois2004}, but Pareto optimality has been considered only scarcely \citep{Warmflash2012}. One issue is that the diversity of environmental conditions is rarely known and the number of potential environments is very large. For example, an \textit{Escherichia coli} bacterium harbours more than 400 regulatory genes responding to typically as many signals. Following an inference approach and testing all potential environments given a fixed gene regulatory architecture would require the consideration of combinations of these signals, the number of which is very large.

In cells, external signals are processed by signal molecules or gene products (e.g. transcription factors) modulating the production or the activity of others. These modulations follow various connectivity patterns comprising cascades and logical integration. While it is not yet possible to fully predict the response of an arbitrary gene network based only on its connectivity, simplified approaches allow classification of the behaviour of networks, among which is the theory of monotone systems \citep{Angeli2003}. A system is said to be monotone if the relation between any pair of input and output is monotone and this monotonicity is independent of the state of the rest of the system.

Gene networks are most often represented by signed graphs, where the sign of the arrow connecting two genes represents a monotone relation between the upstream and downstream gene. It is not guaranteed that the relation between any two arbitrary genes taken within a larger signed network is monotone. Nevertheless, it has been shown that gene networks are essentially monotone both in practice and in theory  \citep{Gjuvsland2013,Sontag2005}, in the sense that: (i) the whole network can be decomposed into a few large monotone components and (ii) the network can be made monotone by removing a small number of genes (a few among hundreds to thousands). It is important to note that the monotonicity described here is established only by the sign of the interactions. Thus, mutations changing the strength of interactions but conserving the signed regulatory architecture leave the monotonicity properties invariant. 

\begin{figure}[ht]
\includegraphics[width=\textwidth]{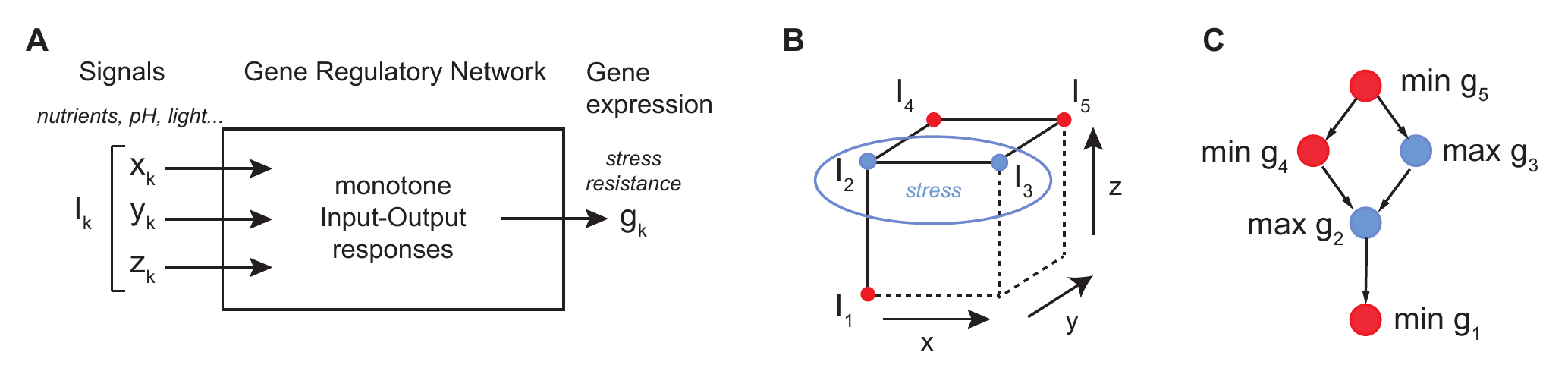}\caption{Evolution of gene regulatory networks in variable environments. \textbf{A)} $I_k$ is a vector of signals $x_k, y_k, z_k, ...$ which are sensed by the organism in environment $k$. Other signals, for example the presence of an antibiotic, may not be sensed. The Gene Regulatory Network is assumed to have monotonous input-ouput relations. Consider the expression of a single gene $g_k$ in the different environments. To illustrate our problem, we assume that g is an antibiotic resistance gene, which provides a fitness benefit in the presence of antibiotics and a cost otherwise (e.g. burden of unnecessary protein expression). \textbf{B)} The variable environment consists of several combinations of sensed signals $I_k$, some of which correspond to environments with antibiotics (in blue). The arrows indicate increasing values $x, y, z$. \textbf{C)}. The natural partial order between coordinates $(x,y,z)\in \mathbb{R}^{3}$ induces a partial order between the $g_k$ values due to the the monotonicity of the response, where, additionally, $g_2$ and $g_3$ must be maximized to resist antibiotics, while the other $g_k$ values must be minimized to limit cost of expression. \label{fig:biology}}
\end{figure}

In this work, we focus our attention on the expression level $g\in \mathbb{R}$ of a given gene in response to a vector of input signals $I \in \mathbb{R}^{k}$ (Figure \ref{fig:biology}). We assume that the response $\mathcal{F}$ between the multiple inputs $I$ and the single output $g$ is monotone. A variable environment can be represented by a list of vectors of input signals $I_1,...,I_N\in \mathbb{R}^{k}$ and their corresponding responses $g_1,...,g_N\in \mathbb{R}$, indexed by  environment. The resulting vector $G = (g_1,...,g_N)\in \mathbb{R}^{N}$ is the vector of expression levels of a single gene in the $N$ environments, $\mathbb{R}^{N}$ being the decision space.

The natural order on $\mathbb{R}$ induces a partial order between the elements $I_i \in \mathbb{R}^{k}$, which, given the monotonicity of $\mathcal{F}$, induces partial order constraints between the values $g_1,...,g_N$. These partial order constraints in turn determine the feasible space for the vector $G\in \mathbb{R}^{N}$. Here we are interested in predicting the evolutionary potential of a fixed signed regulatory network topology. Under the latter condition, mutations alter the strength of the connections, leading to the modification of the response function from $\mathcal{F}$ to $\mathcal{F}'$, but the monotonicity of  $\mathcal{F}$ and $\mathcal{F}'$ is the same. This way, \textit{via} a transformation from $\mathcal{F}$ to $\mathcal{F}'$, mutations occuring during biological evolution move the system from $G$ to $G'$ in the feasible space constained by the partial order. 

We now introduce the optimization objectives: certain coordinates $g_i$ must be maximized while others must be minimized, in accordance with the idea that genes can be considered either detrimental (cost of gene expression) or beneficial, depending on the environment. Searching for Pareto optimal solutions given the partial order constraints finally reveals the set of expression levels $G$ that can be reached during the evolution of a fixed signed monotone regulatory network.

\section{Notations and formulation of the problem}

\subsection{Pareto optimality}

We consider a partially ordered set ($\Omega,\succeq$) with the corresponding strict order $x\succ y\Longleftrightarrow(x\succeq y$ and $x\neq y)$. 
We denote $Par(X)$ the strong Pareto optimal set of $X\subseteq\Omega$.
$Par(X)$ contains the elements of $X$ which are not stricly dominated, in the sense of the partial order above, by another element of $X$:%
\begin{equation}
x\in Par(X)\ \Longleftrightarrow\ x \in X \textrm{and } \forall y \in X \textrm{with } y\succeq x\ \textrm{we have } y=x. \label{eq:maxequal}%
\end{equation}

We note here for our later proof two elementary properties of Pareto fronts, which are general properties of upper sets \citep{Davey2002}: 
\begin{property}
\label{union} 
$\forall A_{1},..., A_{I}  \subseteq\Omega$, $Par\left(\bigcup_{i=1}^{I}A_{i}\right)  \;\subseteq\;\bigcup_{i=1}^{I}Par(A_{i}).$
\end{property}

\begin{property}
\label{necessary} Consider $A, B\subseteq\Omega$, $B$ such that every $y \in B$ is dominated by an element $x \in Par(B)$, and $Par(A) \subseteq B$. Then $Par(A)=Par(A\cap B)$.
\end{property}

The latter property is proved as follows. $Par(A) \subseteq B$ implies $Par(A) \subseteq Par(A \cap B)$. The reverse inclusion comes from using $X=A \cap B$ and $ Y=A$ in: if $X \subseteq Y$ and $Par(Y)\subseteq Par(X)$, then $Par(X)=Par(Y)$. 

For the remainder of this work, we will consider $N$ bounded variables $(x_1, ..., x_N) \in [0,1]^{N}$ and use the convenience notation $x_0 = 0$ and $x_{N+1} = 1$. 

\subsection{Optimization objectives}

Our multi-objective optimization problem consists of minimizing some of the variables and maximizing all the others. Formally, we have a partition $\left\{  \mathcal{N}_{+},\mathcal{N}_{-}\right\}  $ of
the index set $\mathcal{N}=\left\{  1,2,\ldots,N\right\}$. We define a Pareto ordering
$\succeq$ on $\mathbb{R}^{N}$ with the \emph{signature } $\left\{
\mathcal{N}_{+},\mathcal{N}_{-}\right\}  $ as follows: if $x=\left(
x_{1},x_{2},\ldots,x_{N}\right)  \in$ $\mathbb{R}^{N}$and $y=\left(
y_{1},y_{2},\ldots,y_{N}\right)  \in$ $\mathbb{R}^{N}$ then%
\begin{equation}
x\succeq y\Longleftrightarrow\forall i\in\mathcal{N}_{+}:x_{i}\geq y_{i} \text{ and } \forall j\in\mathcal{N}_{-}:x_{j}\leq y_{j},
\end{equation}
i.e. for the variables in the \emph{ascending} set $\left\{  x_{i}%
|i\in\mathcal{N}_{+}\right\}  $ \textquotedblleft larger is
better\textquotedblright, while for those in the \emph{descending} set
$\left\{  x_{j}|j\in\mathcal{N}_{-}\right\}  $ \textquotedblleft smaller is
better\textquotedblright. In other terms, the ideal point has coordinates $0$ for indexes in $\mathcal{N}_{-}$ and $1$ for indexes in $\mathcal{N}_{+}$. We define the corresponding (weak) strict
Pareto order through $x\succ y\Longleftrightarrow x\succeq y\wedge x\neq y$.
Partial order constraints correspond to a set of (weak) inequalities of the type $x_{i}\geq x_{j}$ between coordinates of $x\in\Omega$. 

We consider a sub-class of linear MOPs:

\begin{equation}
\begin{split}
&\textrm{max } x_i, \quad \forall i \in \mathcal{N}_{+}\\
&\textrm{min } x_i, \quad \forall i \in \mathcal{N}_{-}\\
&\qquad \textrm{s.t. } 0 \leq x_i \leq 1, \forall i \in \mathcal{N}\\
&\qquad \textrm{and } x_i \leq x_j, \textrm{for some } i,j \in \mathcal{N}
\end{split}
\end{equation}

This problem has the particularity that the number of objectives equals the number of decision variables. The dimension of the objective space can thus become very large compared to what has typically been considered so far for linear MOPs. We present in the following an algorithm which solves this subclass of large size linear MOPs efficiently.

\subsection{Graph representation of the problem}
\label{definitions}

\noindent\textit{\textbf{Vertex colouring}} --- We introduce a graph stucture whose vertices $V_I$ are associated with a group of variables $\{x_i\}_{i \in I}, I \subseteq \{0, ..., N+1\}$, called \textit{aggregate}, and denoted as $x_I$. An aggregate $x_{I}$ represents a state such that $\forall i,j\in I,x_{i}=x_{j}$, a value noted $x_I$. A vertex can be of one of the four following natures (or colours indicated in parentheses, as used in Figure \ref{fig:illustration}): 
\begin{itemize}
\item if $I \subseteq \mathcal{N}_{-}$, $V_I$ is a \textit{descending vertex} (red);
\item if $I \subseteq \mathcal{N}_{+}$, $V_I$ is an \textit{ascending vertex} (blue); 
\item if $I$ comprises indexes from both $\mathcal{N}_{-}$ and $\mathcal{N}_{+}$, $V_I$ is a \textit{Trade-Off Vertex} or \textit{TOV} (grey); 
\item if $0 \in I$ or $N+1 \in I$, $V_I$ is a \textit{boundary vertex} (black). 
\end{itemize}

\noindent\textit{\textbf{Partial order constraints}} --- Edges express order constraints between variables, according to the convention of so-called Hasse diagrams: edge $E\left(i \rightarrow j\right)$ points from $V_i$ to $V_j$ if and only if $x_{i}\geq x_{j}$. The inequality constraints can be consistently carried over to aggregates provided their index-sets $I$ and $J$ are disjoint, and we will use the obvious notation $x_{I}\geq x_{J}$. Variables engaged in cyclical inequalities are trivially equal to each other and  are thus assumed to be aggregated, resulting in acyclic diagrams.

\noindent\textit{\textbf{Relations between vertices}} --- Two vertices $V_i$ and $V_j$ connected by $E\left(i \rightarrow j\right)$ are said to be \emph{conflicting} if $V_i$ is descending and $V_j$ is ascending. We say that a vertex $V_i$ \emph{aims at} another vertex $V_j$ if: either $V_i$ is descending and points to $V_j$ \textit{via} an edge $E\left(i \rightarrow j\right)$; or vice versa, $V_i$ is ascending and is pointed from $V_j$ \textit{via} an edge $E\left(j \rightarrow i\right)$, independently of the colour of $V_j$. Note that conflicting vertices necessarily aim at each other. A maximal connected component of the Hasse diagram exclusively comprising ascending or descending vertices is called a \textit{monotone connected component}. A vertex is qualified as \emph{extremal} (in the sense of the monotone connected components) if it only aims at vertices of different colour. Note that conflicting vertices are not necessarily extremal as they may point to other vertices of the same colour, and extremal vertices are not necessarily conflicting as they may point to TOVs or to boundary vertices.

\noindent\textit{\textbf{Edge contraction rules}} --- An edge contraction between two vertices $V_I$ and $V_J$ consists of removing $E\left(I \rightarrow J\right)$, and replacing $V_I$ and $V_J$ by a unique vertex $V_K$, where $K = I \cup J$. While all vertices that are not affected by the contraction are of constant colour, the colour of $V_K$ is determined consistently with the colour definitions:
\begin{itemize}
\item if $V_I$ and $V_J$ are ascending, $V_K$ is ascending;
\item if $V_I$ and $V_J$ are descending, $V_K$ is descending;
\item if $V_I$ or $V_J$ is a boundary vertex, $V_K$ is a boundary vertex;
\item if $V_I$ and $V_J$ are different and none of them is a boundary vertex, $V_K$ is a TOV;
\end{itemize}

\noindent\textit{\textbf{Resolution tree}} --- We have just introduced all of the notions for representing our initial problem and any state resulting from edge contractions in the form of 4-colour directed acyclic graphs, for which we reserve the term \emph{Hasse diagrams}, or \emph{diagrams}, denoted $\mathcal{H}$ with an index. The steps of the dynamic algorithm described below generate another type of graph, which we call the \textit{resolution tree}. To avoid any confusion between the two kinds of graphs, the vertices of the resolution tree are called \emph{nodes}. Each node of the resolution tree corresponds to a Hasse diagram, and each edge of the resolution tree corresponds to an operation of edge contraction applied to a Hasse diagram (Figure \ref{fig:illustration}). A \emph{branch} will only refer to a branch of the resolution tree. A vertex aiming at several other vertices in $\mathcal{H}$ will be called a \emph{junction}, whereas a node connected to several downstream nodes in the resolution tree will be called a \emph{branching}. 

\noindent\textit{\textbf{Graph formulation of the optimization problem}} --- Let $\Omega_{\mathcal{H}}\subseteq [0,1]^{N}$ be the space of all vectors respecting the partial order constraints represented by $\mathcal{H}$. Determine $Par(\Omega_{\mathcal{H}})$, the set of vectors optimal under the Pareto order $\succeq$ on $\mathbb{R}^{N}$.

\section{Resolution by edge-contractions of the Hasse diagram}

\subsection{Graph contraction algorithm}

\begin{figure}[ht]
\includegraphics[width=\textwidth]{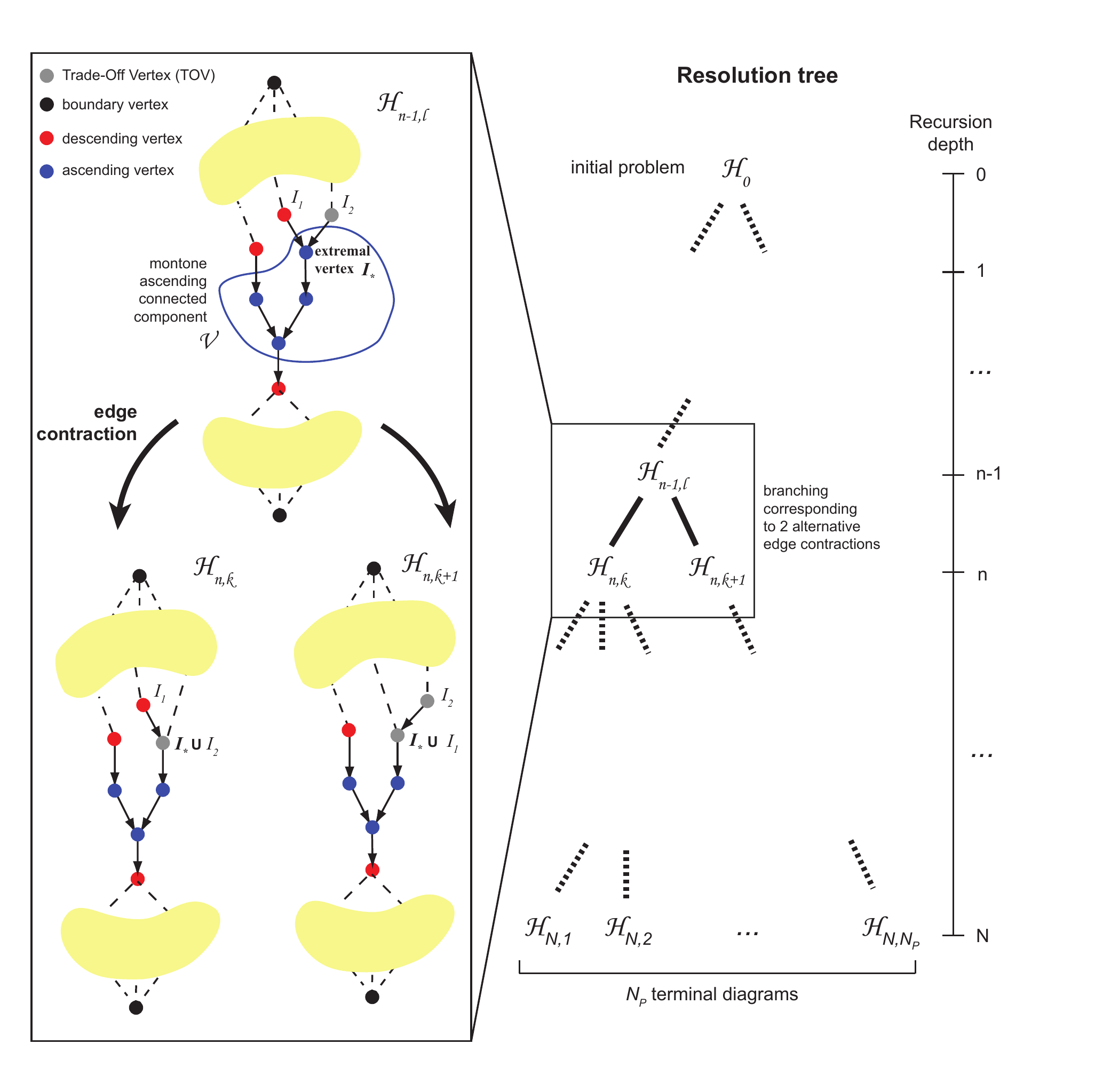}\caption{Graphical illustration of the key steps of our algorithm. \textbf{Left}: Selection of an ascending connected component $\mathcal{V}$ in a Hasse diagram $\mathcal{H}_{n-1,l}$ and identification of an extremal vertex $V_{I_{\ast}}$. $\mathcal{H}_{n-1,l}$ leads to two novel Hasse diagrams $\mathcal{H}_{n,k}$ and $\mathcal{H}_{n,k+1}$ by alternative edge contractions leading to the aggregation of $V_{I_{\ast}}$ with the respectively conflicting vertices $V_{I_1}$ and $V_{I_2}$. \textbf{Right}: The nodes of the resolution tree represent the Hasse diagrams obtained recursively by alternative contractions of single edges (as described on the left). Alternative contractions correspond to branching operations within the resolution tree. The process starts with the initial problem represented by $\mathcal{H}_{0}$ and ends with the terminal diagrams at depth $N$ of the recursion. Terminal diagrams only contain Trade-Off Vertices (TOVs) and variables aggregated to the lower or the upper bound. In the improved version of the algorithm, the terminal diagrams each represent one of the distinct $N_P$ faces of the Pareto front.\label{fig:illustration}}
\end{figure}

The algorithm starts by setting $\mathcal{H}_{0}=\mathcal{H}.$ Steps 1 to 4 described below (and illustrated in Figure \ref{fig:illustration}) are then recursively applied to all diagrams $\mathcal{H}_{n,t}$, $t=1, ..., T_n$, generated at depth $n$ until step 2 can no longer be performed, \textit{i.e.} the diagram in question no longer contains any ascending or descending vertex (equivalently, only contains TOVs and boundary vertices):

\begin{enumerate}
\item Perform a \emph{transitive reduction} of $\mathcal{H}_{n}$,
\textit{i.e.} remove any direct edge $E(u \rightarrow v)$ if there
exists a longer \emph{path} from $V_u$ to $V_v$ on $\mathcal{H}_{n}$.
\item Select an extremal vertex $V_i$.
\item Consider the vertices $V_k,k=1,...,K$ which $V_i$ aims at. There is always at least one such vertex, in limiting cases provided by boundary vertices. Define the diagrams $\mathcal{H}_{n,k},k=1,...,K$ by respectively contracting the edge connecting $V_i$ and $V_k$ according to the colouring rules defined in section \ref{definitions}.
\end{enumerate}

At the end of this branching process, we are left with a collection of
terminal graphs $\mathcal{H}_{t},t=1,...,T$, and we posit that the solution of
the initial problem is:%

\[
Par(\Omega_{\mathcal{H}})=\bigcup_{t=1}^{T}\Omega_{\mathcal{H}_{t}}.
\]

The proof of the algorithm is provided in \ref{sec:proof}.

\subsection{Improved algorithm}
\label{improved}

The algorithm described so far has potential redundancies: (i) duplicated diagrams representing the same parameterization; (ii) diagrams $\mathcal{H}_{i}$ which aggregate the initial variables into a sub-partition of another diagram $\mathcal{H}_{j}$. Case (i) happens when an extremal vertex conflicts with two or more other vertices, because the two corresponding contractions occur in a certain order along a resolution branch, and in another order along another branch. Case (ii) happens when a vertex $V_0$ aims at a TOV $V_1$ and a conflicting vertex $V_2$: along a first resolution branch, $V_0$ aggregates with $V_1$, then $V_2$ aggregates with the resulting TOV, whereas along a second branch, $V_0$ aggregates with $V_2$, resulting in two distinct TOVs. Consequently, the parameterization given by two distinct TOVs in the second branch includes the solution obtained in the first.

A way to fix the redundancies due to case (i) would be to store known nodes of the resolution tree in a hash table \citep{Kambhampati2000}. We propose instead an improved version of the contraction rules, which ensures that every terminal diagram represents a distinct face of the Pareto front and which also removes sub-representations of the Pareto front due to case (ii). 

We define \textit{frozen edges} as edges which cannot be contracted. Furthermore, we impose this property to be inherited downstream of the resolution tree, i.e. a frozen edge remains frozen after contraction of other edges. Otherwise, an edge is qualified as \textit{free}. The improved version of the algorithm consists of modifying step 3 of the original algorithm as follows:\newline

\noindent\textit{\textbf{Modified contraction rule}}: If possible, contract extremal vertices $V_i$ which aim at other vertices \textit{via} a single free edge. Otherwise: (i) contract in priority $V_i$ with conflicting vertices, then with TOVs, and, (ii) for each $k = 2, ..., K$, the edges contracted to obtain $\mathcal{H}_{n,i}$, $i=1,...,k-1$, are frozen in $\mathcal{H}_{n,k}$.\\

As shown in the proof (\ref{sec:proofimproved}), the priority and edge freezing rules ensure that the terminal diagrams $\mathcal{H}_{N,t}, t=1,...,N_P$ describe all distinct faces of the Pareto front. A pseudo-code of the improved algorithm is provided in \ref{sec:pseudocode} and its Python implementation can be downloaded at https://hal.inria.fr/hal-01760120.

\section{Complexity of the algorithm}
\label{discussComplexity}

\subsection{Comparison of running times}

Here, we discuss worst case estimates and test numerically the complexity of different Multi-objective Linear Programs (MOLPs). We call $N$ the number of variables and $N_P$ the number of faces of the Pareto front. We remind the reader that a \textit{face} is defined as a maximal convex subset of the Pareto front \citep{Ecker1980}, which  is itself a subset of the convex polytope $\Omega_{\mathcal{H}}$. 

The running time for improved versions of the Benson algorithm \citep{Ehrgott2012} has been shown to increase exponentially with either the number of faces or vertices of the Pareto front, for a fixed dimension of the decision space \citep{Bokler2017}. In our linear MOP, the number of objectives increases linearly with the number of variables $N$. The running time of Benson's algorithm is thus expected to scale as $O(N_P^N)$. 

For the graph contraction algorithm, each branching operation requires to copy a Hasse diagram of size at most $N$, then to visit at most every vertex to find an extremal one. The diagram can be maintained transitively reduced by only testing for indirect paths connecting the fused vertex and the vertices that were aiming at it before contraction, which is performed in $O(N+E)$ running time with $E$ the number of edges. An upper bound is then $O(N^2)$ operations per branching. From each node of the tree, there are at most $N$ branches. Finally, the size of the resolution tree is bounded by $O(N \cdot N_P)$. Indeed, the tree has a depth of $N$ and each leaf of the tree corresponds to one of the $N_P$ faces of the Pareto front. The resulting worst case running time then scales as $O(N^4\cdot N_P)$.

\begin{figure}
\includegraphics[width=\textwidth]{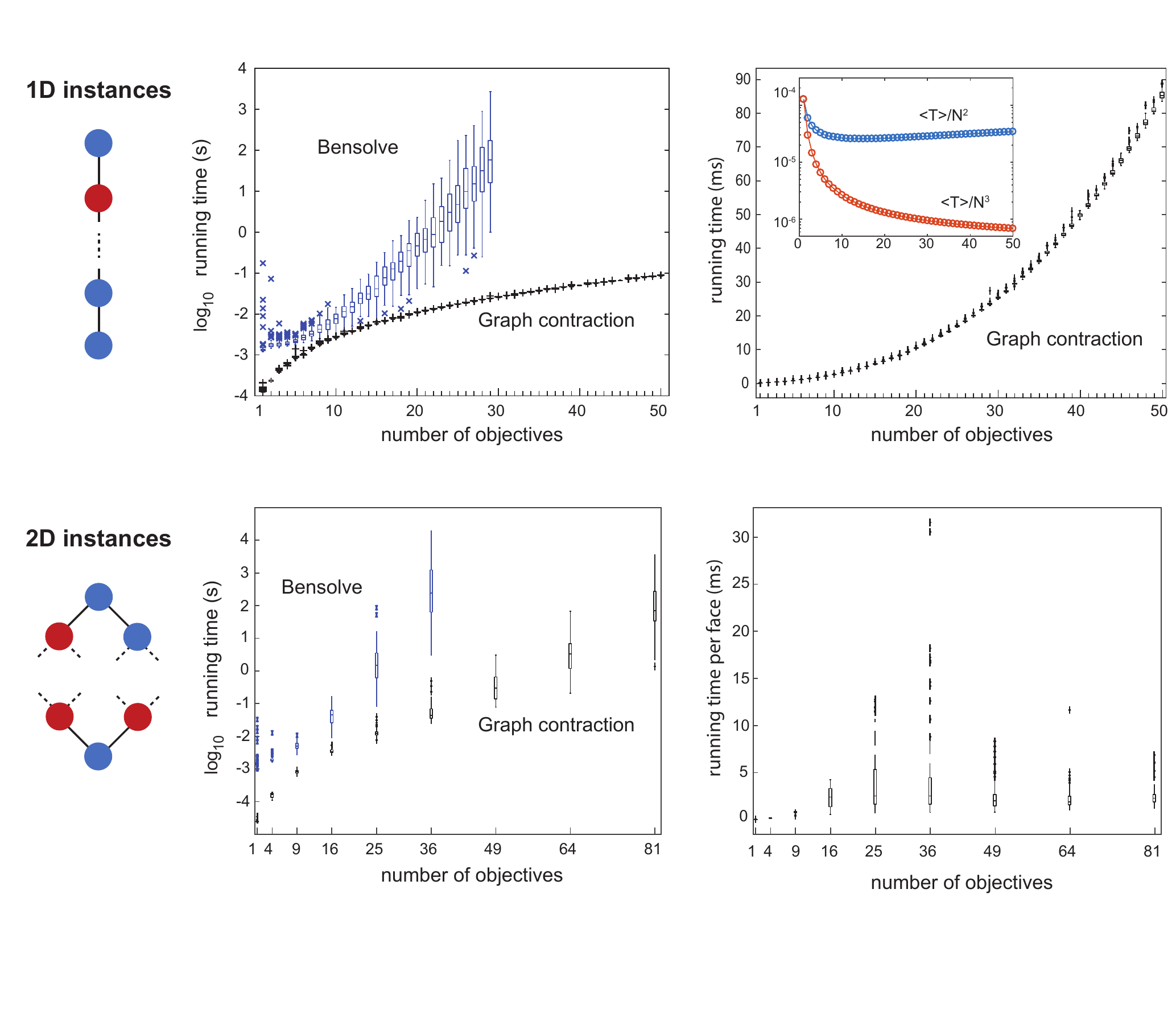}\caption{
Bensolve and the graph contraction algorithm are run on the same sets of 100 random instances for each size, where minimization and maximization objectives have been attributed randomly with probability 0.5. 
\textbf{Top left}: Logarithm of the running times for Bensolve and the graph contraction algorithm in the case of 1D gradients with random stress patches (total order constraints). 
\textbf{Top right}: Same running times on a linear scale for the graph contraction algorithm. Inset: average running time $<T>$ per number of objective divided by $N^2$ (blue, asymptotically increasing) or $N^3$ (red, asymptotically decreasing).
\textbf{Bottom left}: Same as top left graph, but for partially ordered expression levels as imposed by 2D gradients, where the Hasse diagrams correspond to square grids. 
\textbf{Bottom righ}t: Linear scale of the running time per face of the Pareto front for the graph contraction algorithm.
}%
\label{fig:runningtime}%
\end{figure}

We have measured the running time using Bensolve implemented in Matlab and the graph contraction algorithm implemented in Python for the problem of optimization of gene expression in 1D and 2D discretized gradients of environmental signals with random stress patches. 
The corresponding Hasse diagrams are respectively 1D chains (total order,  Figure \ref{fig:runningtime}, top panels) or 2D square lattices (partial order,  Figure \ref{fig:runningtime}, bottom panels) comprising up to 81 objectives. Minimization or maximization objectives are attributed randomly with even probabilities. We tested the same random set of 100 instances for each size of the problem with both algorithms (computer Intel core i7 with 8Go 1600MHz RAM on Windows 10). A table of the data corresponding to Figure \ref{fig:runningtime} are reported in the Supplementary Material.

For the 1D instances, there is no duplication of the Hasse diagrams ($N_P = 1$). As predicted, the running time increases over-exponentially with Bensolve and polynomially for the graph contraction algorithm (Figure \ref{fig:runningtime}, top panels). We estimate the largest polynomial exponent to be between $2$ and $3$ in this case (Figure \ref{fig:runningtime}, top right panel, inset). Similar trends are observed for 2D instances (Figure \ref{fig:runningtime}, bottom left panel). Interestingly, the running time per Pareto face is observed to be almost constant with the graph contraction algorithm (Figure \ref{fig:runningtime}, bottom right panel). This may be due to partial orders on a 2D lattice with randomized objectives leading statistically to monochromatic domains of bounded size as described in the next section.

\subsection{Graph interface and parameterized computation of the Pareto front}

In particular cases of our problem, the graph contraction algorithm can lead to further decrease in computational complexity by using a combinatorial description of the Pareto front.

For this, we introduce a specific subset of the diagram which we call the interface $\mathcal{I}$. This interface comprises all the conflicting vertices of the full diagram. For each resolution of $\mathcal{I}$, the solutions of the monotone connected components can be computed independently, then assembled combinatorially. Parameterization by the different resolutions $\mathcal{I}$ can exponentially reduce the computing time and the size of the description of the full solution. In particular for series-parallel partial orders, the resolution of $\mathcal{I}$ is unique. Under the additional assumption that the size of monotone (or monochromatic) connected components is bounded, one obtains a resolution and a description of the Pareto front in $O(N)$, even though the Pareto front may comprise an exponential number of faces.

We define the interface $\mathcal{I}$ of the initial problem $\mathcal{H}_{0}$ as the set of all conflicting vertices. $\mathcal{I}$ contains all extremal vertices of $\mathcal{H}_{0}$ which do not directly aim at the maximum or minumum bound. While $\mathcal{I}$ can be composed of several connected components, a monotone connected component of $\mathcal{H}_{0}$ may intersect several connected components of $\mathcal{I}$. We call $\mathcal{H}_{u}$ the diagrams obtained by aggregating first all the extremal vertices of $\mathcal{H}_0$ in $\mathcal{I}$. As all conflicting vertices have been aggregated into TOVs at this point, the algorithm only results in aggregation of extremal vertices with existing TOVs. In this sense, the remaining monotone connected components of each $\mathcal{H}_{u}$ are isolated from each other by TOVs. 

Now call $\mathcal{C}_{u,v}, v=1,..,V_u$ the montone connected components of $\mathcal{H}_{u}$ taken together with the TOVs they aim at. Each $\mathcal{C}_{u,v}$ can be solved separately, leading to its own set of leaves $\mathcal{C}_{u,v}^w, w=1,...,W_{u,v}$. The parametrizations of the different parts of the Pareto front of $\mathcal{H}_{u}$ can be obtained by concatenating all possible combinations of the $w$ indexes of the $\mathcal{C}_{u,v}^w$. Here the concatenation $\oplus$ between diagrams is defined as the merging of vertices which represent aggregated variable sets with a non-empty intersection. With this notation, we have:
\begin{equation}
Par(\mathcal{H}) = \bigcup_{u=1}^{U} Par(\mathcal{H}_u)
\end{equation}
with for every $u$:
\begin{equation}
Par(\mathcal{H}_u) = \bigcup_{(w_1, ..., w_v) \in W^\star} \left( \mathcal{C}_{u,1}^{w_1} \oplus ... \oplus \mathcal{C}_{u,v}^{w_v} \right)
\end{equation}
where $W^\star = \{1, ..., W_{u,1} \} \times ... \times \{1, ..., W_{u,v} \}$.

Such a combinatorial representation of the Pareto front can be exponentially smaller than the number of faces of the Pareto front itself. In particular when the size of the $\mathcal{C}_{u,v}$ is bounded, the number of terms of the concatenation representing $Par(\mathcal{H}_u)$ increases linearly, while they represent an exponentially increasing number of faces of the Pareto front. 

In the case of series-parallel diagrams $\mathcal{H}$, we have $U=1$. This is due to the forbidden sub-graph characterization of series-parallel graphs: fence subgraphs (``N" shaped motifs) are absent. This property implies that there cannot be conflicting vertices which each participate to a junction. In other words, at least one of the two has no other alternative than contracting with its conflicting vertex, leading to the absence of branching process during the resolution of the interface of the Hasse diagram. Under the additional condition that monochromatic connected components have bounded size, all $\mathcal{C}_{1,v}$ have a bounded size, and the complexity of the resolution of these sub-diagrams is bounded. Thus the Pareto front of the full problem admits a representation which complexity grows linearly with the number of $\mathcal{C}_{1,v}$, which itself increases at most linearly with the number of initial variables.

\section{Applications}

\subsection{Example of stress patches in 2D gradients}

\begin{figure}
\includegraphics[width=\textwidth]{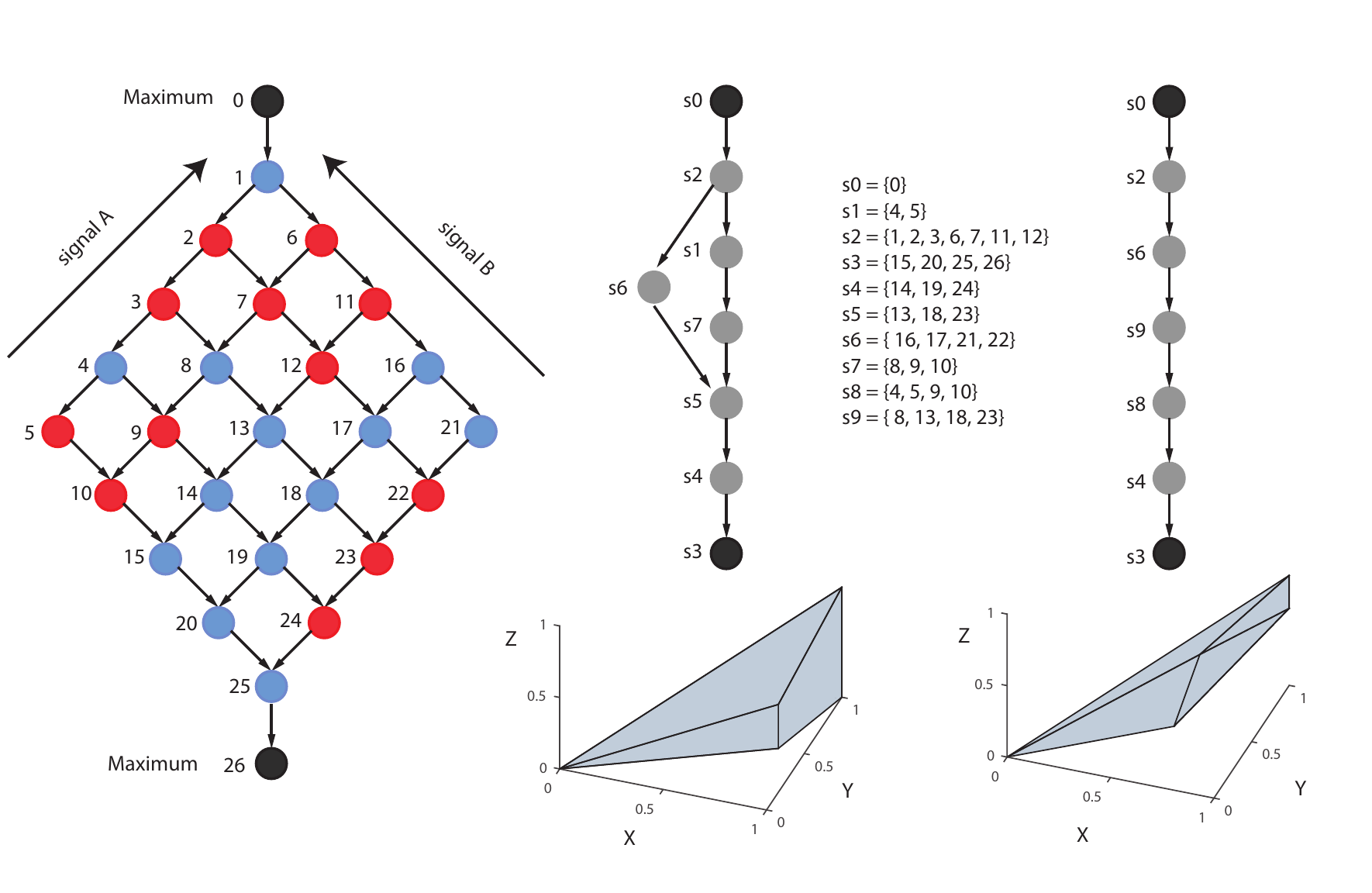}\caption{
Example of resolution of Pareto optimal gene expression for a monotonous system in response to combinations of 5 levels signals A and B, with three stress patches (red vertices corresponding to gene expression maximization objectives). On the right are represented the Hasse diagrams of the two Pareto faces, on top of their respective projections along the unit vectors X, Y and Z,  corresponding to the directions defined by $x_i = 1$ for $i = 1, ..., 8$, $0$ otherwise,  $y_i = 1$ for $i = 9, ..., 12$, $0$ otherwise,  $z_i = 1$ for $i = 13, ..., 25$, $0$ otherwise. Indexed sets s1 to s9 are aggregated variables, which are subsets of variables equal to each other within a face of the Pareto front.
}%
\label{fig:example}%
\end{figure}

To illustrate our algorithm, we show an example of Pareto optimal solutions in the case of a 2D gradient (Fig. \ref{fig:example}). In this problem, all combinations of five levels of each signal (A and B) result in 25 environments. We consider a monotonous response regulatory network, leading to partial order constaints between expression levels that reflect the partial order between signal intensities. Gene expression must be minimized in the absence of stress (cost of gene expression) and maximized in the stress patches (red vertices in Fig. \ref{fig:example}, benefit of gene expression, e.g. antibiotic resistance gene in the presence of antibiotics). Computation shows that the Pareto optimal set comprises two distinct faces of dimensions 5 and 6.

\subsection{Exhaustive characterization of Pareto fronts}

\begin{figure}
\includegraphics[width=\textwidth]{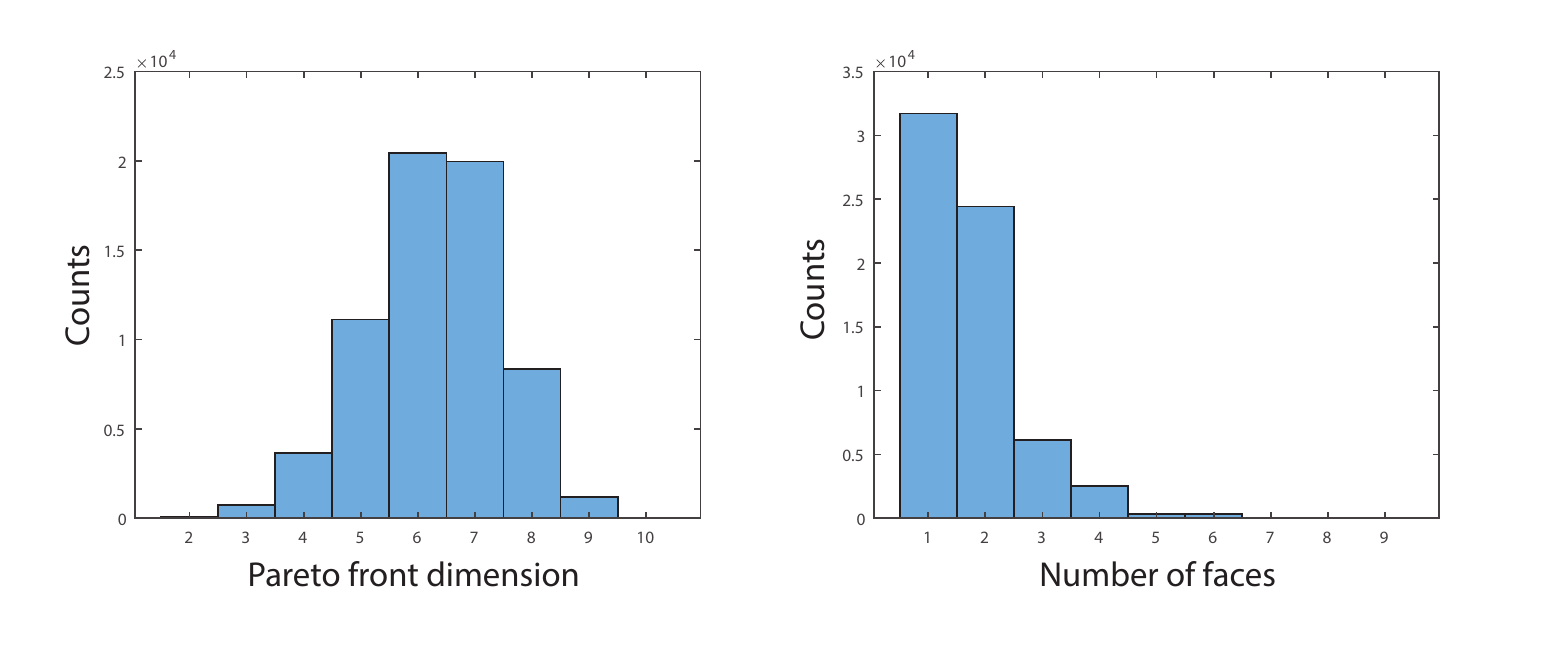}\caption{
Exhaustive characterization of optimal solutions in 2D gradients with 4 signal intensities leading to $N=16$ environments. The dimension of the Pareto front (left hand panel) and the number of faces (right hand panel) are computed for all possible combinations of maximization and minimization objectives.
}%
\label{fig:screen}%
\end{figure}

The Pareto optimal sets corresponding to all stress patches configurations in 2D gradients comprising $N=16$ environments (all combinations between 4 levels of 2 signals) have been computed. Figure \ref{fig:screen} shows the distribution of Pareto front dimensions (maximum face dimension) and of the number of faces. The corresponding bivariate table is provided in the Supplementary Material.

\section{Discussion and conclusion}

We have introduced a multi-objective optimization problem relevant to the evolution of gene regulatory networks, where the monotone nature of the input-output response leads to partial order constraints between the level of expression of a gene in different environments. In turn, selective pressures tend to maximize or minimize the expression of the gene in the variable environment, such that the number of objectives is equal to the number of different environments. More generally, this approach can provide bounds on the optimal solutions accessible by monotone regulatory systems.

We have described and demonstrated an algorithm which provides an exact parameterization of the full Pareto front for this problem. This corresponds geometrically to finding a polytope defined by partial order relations within the hypercube $[0,1]^N$ given an ideal point located at a corner of this hypercube. The solution is obtained in polynomial time with the number of objectives $N$ (which equals the number of variables in our specific problem) relative to the number $N_P$ of faces of the Pareto front. The graph contraction algorithm exploits the particularities of this linear MOP, and scales polynomially with the number of objectives, in contrast to general Multi-Objective Linear Programs (MOLPs), in which complexity scales at least exponentially.

To achieve this, we established a mapping between Hasse diagrams and polytopes, whereby a colouring of the graph encodes the location of the ideal point. In this approach, vertices represent sets of aggregated variables, edges correspond to ordering relations, and colours correspond to the optimization objectives associated with the variables. Following a dynamic programming approach in the space of coloured graphs, the initial polytope is recursively reduced to lower-dimension spaces, corresponding to edge contractions in its diagrammatic representation. The Pareto front ultimately consists of the union of spaces corresponding to the terminal Hasse diagrams obtained after $N$ contractions. 

We have furthermore introduced a parmeterized complexity approach, by considering a specific subgraph, which we call the interface and which corresponds to the smallest set containing all the potential trade-offs. The edge contraction algorithm can be applied to this subgraph prior to the rest of the graph. For each resolution of the interface, the remaining coloured connected components can then be solved independently of one another, and the Pareto front can be represented combinatorially from the solutions of these isolated components. Interestingly, for series-parallel partial orders, the interface has a unique resolution. When, additionally, the remaining coloured connected components are of bounded size, the Pareto front, though of exponential complexity, can be computed and represented in $O(N)$.

A source of efficiency of the graph contraction algorithm is the exploitation of the fact that each constraint is related to only two objectives (even though every objective or variable can be related to an arbitrary number of constraints). It would be interesting to investigate whether the process of conflict resolution used here could be applied more generally to efficiently solve linear MOPs with a large number of objectives but where each constraint is associated with a bounded number of them.

\section*{Acknowledgments}
PN was supported by FOM Programme 103 ``DNA in Action: Physics of the Genome''. This work is part of the research programme of the Foundation for Fundamental Research on Matter (FOM), which is part of the Netherlands Organisation for Scientific Research (NWO). We thank Olivier Spanjaard from the Laboratoire d'Informatique de Paris 6 (LIP6) and the members of the Laboratoire d'Informatique Fondamentale de Lille (LIFL) for discussions, and Harry Kemble for careful reading of the manuscript. The author also thank the anonymous reviewers for their contributions.





\appendix

\section{Proof of the graph contraction algorithm\label{sec:proof}}

This proof proceeds in three steps:

\begin{description}
\item[Step 1] We first show that for each iteration of the algorithm,
$Par(\Omega_{\mathcal{H}_{n}})=Par(\bigcup_{k=1}^{K}\Omega_{\mathcal{H}_{n,k}})$.

\item[Step 2] Next we show that $Par(\bigcup_{k=1}^{K}\Omega_{\mathcal{H}%
_{n,k}})=\bigcup_{k=1}^{K}Par(\Omega_{\mathcal{H}_{n,k}})$.

\item[Step 3] Finally we show that the terminal graphs satisfy $Par(\Omega
_{\mathcal{H}_{t}})=\Omega_{\mathcal{H}_{t}}$.
\end{description}

For brevity's sake we treat only the case of an ascending connected component, the
demonstration being easily adapted for a descending connected component.\newline

\textbf{Step 1}

Consider an ascending subgraph $\mathcal{V}$ in $\mathcal{H}_{n}$ and a
maximal vertex $V(x_{I_{\ast}})\in\mathcal{V}$. There are two possibilities:

\begin{enumerate}
\item The vertex upper bound, $V_{N+1}$, is directly connected to $V_{I_{\ast}}$. 
In this case $V_{N+1}$ is the unique vertex pointing to $V_{I_{\ast}}$. 
Indeed, if there would be another vertex $V$ pointing to $V_{I_{\ast}}$, there would be a chain of vertices pointing from $V_{N+1}$ to $V$ then $V_{I_{\ast}}$ (as $V_{N+1}$ is a global upper bound). 
This would contradict the fact that we have taken the transitive reduction of the diagram.

\item Otherwise, call $V_{I_k},k=1,...,K$ the vertices pointing to $V(x_{I_{\ast}})$.
\end{enumerate}

In both cases, we define the sets
\[
\omega\left(  I_{\ast}|I_k\right)  =\left\{  x\in\mathbb{R}^{N}\ \text{such
that}\ x_{I_{\ast}}=x_{I_k}\right\}  .
\]

A necessary condition for $x\in Par(\Omega_{\mathcal{H}_{n}})$ is $x\in
\bigcup_{k=1}^{K}\omega\left(  I_{\ast}|I_{k}\right)  $. Suppose otherwise
that $\forall k:x_{I_{\ast}}<x_{I_{k}}$ (in the case 1, $x_{I_{\ast}%
}<x_{N+1}$). There would exist $\epsilon>0$ such that for all $k:x_{_{I_{\ast}}%
}+\epsilon<x_{I_{k}}$. If we denote by $\epsilon_{I_{\ast}}\in\mathbb{R}^{N}$
the vector with coordinates $x_{j}=\epsilon$ for $j\in I_{\ast}$ and otherwise
$0$, we have for all $x+\epsilon_{I_{\ast}}\in\Omega_{\mathcal{H}_{n}}$ and
$x+\epsilon_{I_{\ast}}\succ x$ as by assumption $I_{\ast}\subseteq
\mathcal{N}_{+}$. This contradicts $x\in Par(\Omega_{\mathcal{H}_{n}}%
)$.\newline

As by definition, $\Omega_{\mathcal{H}_{n}}\cap\omega\left(
I_{\ast}|I_{k}\right)  =\Omega_{\mathcal{H}_{n,k}}$,  we can use Proposition \ref{necessary} and have that $Par(\Omega
_{\mathcal{H}_{n}})=Par(\Omega_{\mathcal{H}_{n}}\cap\bigcup_{k=1}^{K}%
\omega\left(  I_{\ast}|I_{k}\right)  )=Par(\bigcup_{k=1}^{K}\Omega
_{\mathcal{H}_{n,k}})$. \newline

\textbf{Step 2}

We only have to show that $Par(\bigcup_{k=1}^{K}\Omega_{\mathcal{H}_{n,k}%
})\supset\bigcup_{k=1}^{K}Par(\Omega_{\mathcal{H}_{n,k}})$, the inclusion in
the other direction follows directly from Proposition \ref{union}. We
only discuss the case $k\geq2$ as the result is trivial otherwise (in
particular in case 1 of the first step).\newline

Consider two distinct indexes from the set $\left\{  1,..,K\right\}  ,$ which
we can take to be $1$ and $2$ without loss of generality. Consider $x\in
\Omega_{\mathcal{H}_{n,1}}$ and suppose there exists $y\in\Omega
_{\mathcal{H}_{n},_{2}}$ such that $y\succeq x$. Then by virtue of $V_{I_{\ast}}$
being an ascending vertex and by definition of the Pareto order:
\begin{equation}
y_{I_{\ast}}\geq x_{I_{\ast}} \label{c1}%
\end{equation}
Next, as $V_{I_{\ast}}$ was maximal within its ascending subgraph, any
vertex pointing to it must contain at least one descending variable, from
which it follows that
\begin{equation}
x_{I_{1}}\geq y_{I_{1}}. \label{c2}%
\end{equation}
Indeed, the latter inequality is trivially true if all $i\in I_{1}$ label descending
variables, i.e. $I_{1}\subseteq\mathcal{N}_{-}$. If not, we can choose an
$i_{+}\in I_{1}\cap\mathcal{N}_{+}$ and $i_{-}\in I_{1}\cap\mathcal{N}_{-}$
for which the statement $y\succeq x$ implies both $y_{i_{+}}\geq x_{i_{+}}$
and $x_{i_{-}}\geq y_{i_{-}}$ and hence, together with $x_{i_{+}}=x_{i_{-}}$
and $y_{i_{+}}=y_{i_{-}}$ by definition of the aggregates, we obtain the equality $x_{I_{1}%
}=y_{I_{1}}$.

We also have $x\in\Omega_{\mathcal{H}_{n,1}}\subseteq\omega\left(  I_{\ast
}|I_{1}\right)  =\left\{  z\in\mathbb{R}^{N}\ \text{such that}\ z_{I_{\ast}%
}=z_{I_{1}}\right\}  $ implying:
\begin{equation}
x_{I_{\ast}}=x_{I_{1}}. \label{c3}%
\end{equation}
Finally, $y\in\Omega_{\mathcal{H}_{n,2}}\subseteq\Omega_{\mathcal{H}_{n}}$ and
in $\mathcal{H}_{n}$, vertex $V_{I_{1}}$ points to $V_{I_{\ast}}$ by hypothesis, which implies:
\begin{equation}
y_{I_{1}}\geq y_{I_{\ast}}. \label{c4}%
\end{equation}
Examining the above relations in the order (\ref{c2}-\ref{c4}-\ref{c1}%
-\ref{c3}), we see that all the variables at play must be equal, in particular
$y_{I_{1}}=y_{I_{\ast}}$. Then $y\in\omega\left(  I_{\ast}|I_{1}\right)  $ and
consequently $y\in\Omega_{\mathcal{H}_{n,1}}$. To summarize, we have just
demonstrated:
\begin{equation}
\forall j\in\left\{  1,...,K\right\}  ,\forall x\in\Omega_{\mathcal{H}_{n,j}%
}:(y\in\bigcup_{k=1}^{K}\Omega_{\mathcal{H}_{n,k}}\ \text{and}\ y\succeq
x)\Rightarrow(y\in\Omega_{\mathcal{H}_{n,j}}). \label{inter}%
\end{equation}
Now, if we take in particular $x\in Par(\Omega_{\mathcal{H}_{n,j}})$ in
relation (\ref{inter}), relation (\ref{eq:maxequal}) implies that $y=x$ by
maximality of $x$ in $\Omega_{\mathcal{H}_{n,j}}$. This gives:
\begin{equation}
\forall j\in\left\{  1,...,K\right\}  ,\forall x\in Par(\Omega_{\mathcal{H}%
_{j}}),\forall y\in\bigcup_{k=1}^{K}\Omega_{\mathcal{H}_{k}}:\ (y\succeq
x\Rightarrow y=x).
\end{equation}
Applying this time relation (\ref{eq:maxequal}) in the backward direction demonstrates
the announced result:
\begin{equation}
\forall j\in\left\{  1,...,K\right\}  :Par(\Omega_{\mathcal{H}_{n,j}})\subseteq
Par(\bigcup_{k=1}^{K}\Omega_{\mathcal{H}_{n,k}}).
\end{equation}
\newline

\textbf{Step 3}

By construction a terminal graph $\mathcal{H}_{t}$ contains only boundary vertices and TOVs.
Now consider $x,y\in\Omega_{\mathcal{H}_{t}},y\succeq x$, and consider a
vertex $V_{I}$ in $\mathcal{H}_{t}$.
If $V_{I}$ is a boundary vertex, then the variables in $I$ are already at
their optimum bounds, and $x_{I}=y_{I}.$ Otherwise, $V_I$ is a TOV and
as in Step 2 above we can choose an $i_{+}\in I\cap\mathcal{N}_{+}$ and
$i_{-}\in I\cap\mathcal{N}_{-}$ for which the statement $y\succeq x$ implies
both $y_{i_{+}}\geq x_{i_{+}}$ and $x_{i_{-}}\geq y_{i_{-}}$ and hence
$x_{I}=y_{I}$. 
As the aggregates form a partition of the initial index set $\mathcal{N}$, the above immediately implies $x=y$, and hence, by relation (\ref{eq:maxequal}), $x\in Par(\Omega_{\mathcal{H}_{t}})$. As this
is true for any $x\in\Omega_{\mathcal{H}_{t}}$, we have demonstrated our
result: $\Omega_{\mathcal{H}_{t}}=Par(\Omega_{\mathcal{H}_{t}})$.

\section{Proof of the improved algorithm\label{sec:proofimproved}}

A potential issue with the frozen edge rule would be the creation of an extremal vertex not connected to any free edge, thus prematurely stopping the algorithm. This cannot happen due to the priority contraction of single free edges, which implies that frozen edges are generated only at stages when every extremal vertex is connected to two or more free edges. Therefore, the algorithm can be consistently run until the leaf of each branch is reached.

As the treatment of single free edge contractions does not differ from the rules of the initial algorithm, we set ourselves at a branching of the resolution tree corresponding to $k \geq 2$ alternative contractions of an extremal vertex $V_0$ aiming at $V_1, ..., V_k$, where the $j$ first $V_i$ conflict with $V_0$, and the remaining $V_i$ are TOVs, where $0 \leq j \leq k$. 

We show first that the edge freezing rules lead to all admissible (in the sense of other rules) and distinct partitions of the initial variables. Consider the first contraction of $E\left(0\rightarrow 1\right) $. The resulting graph $\mathcal{H}_{n,1}$ induces all admissible partitions such that $x_0$ and $x_1$ are in the same set. Consider then the contraction of $E\left(0\rightarrow 2\right)$, where $E\left(0\rightarrow 1\right)$ is frozen according to the modified algorithm. The resulting branch $\mathcal{H}_{n,2}$ induces all partitions such that $x_0$ and $x_2$ are in the same set but $x_1$ is not. For each $i$-th iteration of this process, the sub-tree stemming from $\mathcal{H}_{n,i}$ induces partitions such that $x_0$ and $x_i$ are in the same aggregate but $x_1,...,x_{i-1}$ are not. Therefore, throughout the different branches, the contractions with $V_0$ enumerate without redundancy all accessible subsets of $\{x_0,...,x_k\}$ containing $x_0$. 

We now want to show that the space parameterized by every terminal diagram is a distinct face of the Pareto front. We first show that none of the partitions are included in another. This is obtained thanks to the prioritization of conflicting vertices contractions: when $V_0$ aggregates with a TOV $V_l$ with $l > j$, any edge $E(0 \rightarrow i)$ with $i \leq j$ joining $V_0$ to a conflicting $V_i$ is frozen in $\mathcal{H}_{n,l}$. Consequently, at this stage, none of the conflicting $V_i$ can be aggregated to the TOV resulting from the contraction of $E(0 \rightarrow l)$. This implies that two conflicting vertices susceptible to form a TOV cannot both aggregate to another TOV, at any stage of the process. Therefore, a TOV can only contain a single pair of conflicting variables, whereas at least two such pairs would be necessary to form a sub-partition of the aggregate.

Finally, we consider two terminal graphs $\mathcal{H}_1$ and $\mathcal{H}_2$, respectively defining spaces $H_1$ and $H_2$, and a dimension larger than 3, as other cases are trivial. 
The partitions associated with $\mathcal{H}_1$ and $\mathcal{H}_2$ not being included in one another implies the existence of variables $(x,y,z)$ such that: (i) in $\mathcal{H}_1$, $x$ and $y$ are in a same set of the partition while $z$ is in another set; (ii) in $\mathcal{H}_2$, $y$ and $z$ are in a same set of the partition while $x$ is in another set. 
Considering points $A\in H_1 \setminus H_2$ and $B \in H_2 \setminus H_1$, the $(x,y,z)$ coordinates of interior points of the segment joining A and B are strictly different from each other. This shows that $ H_1 \cup H_2$ is non-convex. Terminal graphs thus parameterize maximally convex subsets of the Pareto front.

\section{Peusdo-code \label{sec:pseudocode}}

\textbf{function} \textsc{Pareto} (G)

\textbf{input}: Graph $G$ with list $V_G$ of $N$ vertices

\textbf{output}: list of graphs $OUT$

1:\space\space $i \leftarrow 0$

2:\space\space $condition \leftarrow \textsc{false}$ 

3:\space\space \textbf{repeat}

4:\space\space	\quad $i \leftarrow i+1$
	
5:\space\space	\quad $V \leftarrow V_G(i)$
	
6:\space\space	\quad \textbf{if} $V$ is maximal \textbf{and} $V$ aims at a single vertex \textbf{then}
	
7:\space\space		\quad \quad $V_f \leftarrow V$

8:\space\space		\quad \quad $condition \leftarrow \textsc{true}$
		
9:\space\space	\quad \textbf{end if}
	
10: \textbf{until} $condition$ \textbf{or} $i > N$

11: $i \leftarrow 0$

12: \textbf{while} $\neg condition$ \textbf{and} $i \leq N$ \textbf{do}

13:	\quad $i \leftarrow i+1$
	
14:	\quad $V \leftarrow V_G(i)$
	
15:	\quad \textbf{if} $V$ is maximal \textbf{then}
	
16:		\quad \quad $V_f \leftarrow V$

17:		\quad \quad $condition \leftarrow \textsc{true}$
	
18:	\quad \textbf{end if}
	
19: \textbf{end while}

20: \textbf{if} $\neg condition$ \textbf{then}

21:	\quad $OUT \leftarrow G$
	
22: \textbf{else}

23:	\quad Set $OUT$ to empty list

24:	\quad Set $T$ to empty list 
	
25:	\quad \textbf{for all} $V \in V_G$ \textbf{do}
	
26:		\quad \quad \textbf{if} $V_f$ aims at $V$ \textbf{then}
		
27:			\quad \quad \quad Append $V$ to $T$
			
28:		\quad \quad \textbf{end if}
		
29:	\quad \textbf{end for}
	
30:	\quad Sort $T$ with ascending and descending vertices first
	
31:	\quad \textbf{for all} $V \in T$ \textbf{do}
	
32:		\quad \quad $e \leftarrow$ edge between $V_f$ and $V$ in $G$
		
33:		\quad \quad $G_i \leftarrow$ graph obtained by contracting $e$ in $G$
		
34:		\quad \quad $G_i \leftarrow$ transitive reduction of $G_i$

35:		\quad \quad Append \textsc{Pareto}($G_i$) to $OUT$
		
36:		\quad \quad Freeze $e$ in $G$
		
37:	\quad \textbf{end for}	
	
38: \textbf{end if}

39: \textbf{return} $OUT$

\end{document}